\documentclass[final]{siamltex}

\usepackage{graphicx}
\usepackage{enumitem}

\usepackage{color}
\usepackage{bm}
\usepackage{amsfonts}
\usepackage[tight,footnotesize]{subfigure}

\newcommand{\tr}{^{\sf T}}
\newcommand{\m}[1]{{\bf{#1}}}

\newcommand{\C}[1]{{\cal {#1}}}

\renewcommand{\bar}{\overline}
\renewcommand{\tilde}{\widetilde}

\title{Algorithm XXXX: A Gradient-Based Implementation of the Polyhedral Active Set Algorithm
\thanks{
February 11, 2022.
Revised August 2, 2022.
The authors gratefully acknowledge support by the National
Science Foundation under grants 1819002, 1819161, 2031213, and 2110722, and
by the Office of Naval Research under grants N00014-15-1-2048,
N00014-18-1-2100, and N00014-22-1-2397.
}}
\author{
    William W. Hager\thanks{{\tt hager@ufl.edu},
        http://people.clas.ufl.edu/hager/,
        PO Box 118105,
        Department of Mathematics,
        University of Florida, Gainesville, FL 32611-8105.
        Phone (352) 294-2308.}
\and
    Hongchao Zhang\thanks{{\tt hozhang@math.lsu.edu},
        http://www.math.lsu.edu/$\sim$hozhang,
        Department of Mathematics,
        Louisiana State University, Baton Rouge, LA 70803-4918.
        Phone (225) 578-1982. Fax (225) 578-4276.}
}
\begin{document}
\maketitle
\begin{abstract}
The Polyhedral Active Set Algorithm (PASA) is designed to optimize a general
nonlinear function over a polyhedron.
Phase one of the algorithm is a nonmonotone gradient projection algorithm,
while phase two is an active set algorithm that explores faces of the
constraint polyhedron.
A gradient-based implementation is presented, where a projected version of
the conjugate gradient algorithm is employed in phase two.
Asymptotically, only phase two is performed.
Comparisons are given with IPOPT using polyhedral constrained
problems from CUTEst and the Maros/Meszaros quadratic programming test set.
\end{abstract}


\pagestyle{myheadings}
\thispagestyle{plain}
\markboth{W. W. HAGER AND H. ZHANG}
{GRADIENT-BASED POLYHEDRAL ACTIVE SET ALGORITHM}
\section{Introduction}
The polyhedral active set algorithm PASA is designed to solve the problem
\begin{equation}\label{P}
\min f(\m{x}) \quad \mbox{subject to} \quad \m{x} \in \Omega,
\end{equation}
where $f: \mathbb{R}^n \rightarrow \mathbb{R}$ and $\Omega$ is a polyhedron.
Throughout the paper, it is assumed that
\begin{equation}\label{omega}
\Omega = \{ \m{x} \in \mathbb{R}^n : \m{Ax} \le \m{b} \},
\end{equation}
where $\m{A} \in \mathbb{R}^{m \times n}$ and $\m{b} \in \mathbb{R}^m$.
The PASA software, on the other hand, utilizes the representation
\begin{equation}\label{p_omega}
\Omega = \{ \m{x} \in \mathbb{R}^n :
\m{bl} \le \m{Ax} \le \m{bu}, \quad \m{lo} \le \m{x} \le \m{hi} \},
\end{equation}
with $\m{bl}$ and $\m{bu} \in \mathbb{R}^m$
and $\m{lo}$ and $\m{hi} \in \mathbb{R}^n$;
any of the inequalities could be vacuous.
The software is designed to exploit sparsity in $\m{A}$.

The algorithms implemented in PASA have been developed over more than 20 years.
In one series of papers
\cite{ChenDavisHagerRajamanickam09, DavisHager99, DavisHager01, DavisHager05,
DavisHager09},
Timothy Davis and William Hager developed techniques for
modifying a sparse Cholesky factorization of a matrix of the form
$\m{AA}\tr$ after adding or deleting a small number of columns and rows from
$\m{A}$.
These update/downdate techniques are optimal in the sense that
their running time is proportional to the number of nonzeros in the Cholesky
factorization that change.
In a series of papers
\cite{DavisHager08b, DavisHager08, Hager92, H93, Hager99,
Hager02c, Hager03, HagerHearn93}, Hager developed the
Dual Active Set Algorithm (DASA), first in a general context, and then
with Donald Hearn, it was applied to quadratic network optimization;
later with Timothy Davis \cite{DavisHager08b, DavisHager08}, it was applied
to linear programming using the newly developed update/downdate techniques.
More recently, in \cite{hz16} both DASA and the update/downdate techniques
were used in an algorithm PPROJ to project a point onto a polyhedron.

In another series of papers, \cite{hz03, hzACM04, hz05}
William Hager and Hongchao Zhang developed a fast version
of the conjugate gradient method known as CG\_DESCENT since the
search directions were always descent directions, independent of the line
search.
In \cite{HagerZhang13c}, CG\_DESCENT was enhanced using limited memory
techniques.
As an application of CG\_DESCENT,
an active set method for purely bound constrained problems was developed in
\cite{hz05a, hz05b}.
This active set algorithm had two phases, in phase one the gradient projection
algorithm and a cyclic Barzilai/Borwein \cite{bb88,dhsz05} step
were used to identify active constraints, and in phase two,
an unconstrained solver, such as CG\_DESCENT,
optimized the objective over faces of the polyhedral constraint.
Whenever a new constraint in the polyhedron became active,
the optimization was restricted to the resulting smaller face of the polyhedron.

The polyhedral active set algorithm PASA in \cite{HagerZhang16} is
a generalization of the two phase algorithm in \cite{hz05a} from bound
constraints to polyhedral constraints.
Under nondegeneracy type assumptions,
only the second phase is executed asymptotically;
consequently, the asymptotic convergence speed coincides
with that of the algorithm used to optimize the objective
over the faces of a polyhedron.
For a general polyhedron, the projected gradients of phase one are computed
using PPROJ, while in the special case where the polyhedron is a knapsack-type
constraint
\[
\{ \m{x} \in \mathbb{R}^n : bl \le \m{a}\tr \m{x} \le bu, \;\;
\m{lo} \le \m{x} \le \m{hi} \}, \quad \m{a} \in \mathbb{R}^n,
\]
the projection is computed using the Newton/heap-based algorithm
NAPHEAP of \cite{DavisHagerHungerford16}.
The current version of PASA uses a projected conjugate gradient iteration
in phase two to optimize over a shrinking series of faces of the polyhedron.
%
\section{Literature Review}
\label{literature}
%
We briefly summarize continuous nonlinear optimization algorithm
development during the past 30 years.
Early codes in this timeframe include
MINOS \cite{Minos87},
NPSOL \cite{NPSOL}, OPTPACK \cite{H87, H93}, and LANCELOT \cite{cgt92}.
Murtagh and Saunders' MINOS is based on Robinson's algorithm
\cite{Robinson72} which is locally quadratically convergent.
The Lagrangian in Robinson's algorithm is replaced by an
augmented Lagrangian, and the subproblem associated with the linearized
constraints are solved by a reduced gradient algorithm combined with
a quasi-Newton method as described in
\cite{MurtaghSaunders78,MurtaghSaunders82}.
Gill, Murray, Saunders, and Wright's
NPSOL is a sequential quadratic programming method (SQP) where
a positive definite quasi-Newton approximation to the true
Lagrangian Hessian is utilized.
The resulting quadratic programming problem is solved by codes in
the LSSOL package \cite{LSSOL}, which employs active set methods and
dense linear algebra to solve constrained linear least-squares
problems and convex quadratic programming problems.
OPTPACK \cite{H87, H93}
alternates between a constraint step based on Newton's method
and an optimization step based on the minimization of an
augmented Lagrangian over linearized constraints.
The combined steps are locally quadratically convergent, and
the implementation employs dense linear algebra.
Conn, Gould, and Toint's LANCELOT treats nonlinear constraints
using an augmented Lagrangian
which is minimized within a region defined by the bound constraints.
The bound constrained problem is solved by an algorithm that combines
projected gradient techniques \cite{CGT88}
and special structures to exploit the
group partially separable structure of a problem \cite{CGT90}.

More recently, NPSOL was the starting point for SNOPT
(Sparse Nonlinear Optimizer) \cite{gms05} and DNOPT (Dense Nonlinear Optimizer)
\cite{gsw16}.
Again, a positive definite quasi-Newton approximation to the true
Lagrangian Hessian is employed.
However, in SNOPT sparse linear algebra is used to solve the resulting
quadratic program, while DNOPT employs dense linear algebra.
A different SQP algorithm is developed by Fletcher and Leyffer
in \cite{fl02, flt02},
where a trust region approach is applied to the quadratic programming
problem, and the accepted iterates are chosen using a filter method.
The authors of the augmented Lagrangian-based code LANCELOT changed their
focus to efficient quadratic programming (QP) solvers that could be used in the
implementation of SQP methods.
The new software formed the package GALAHAD.
Further development of reliable augmented Lagrangian techniques for general
nonlinear optimization were continued by Birgin, Mart\'{\i}nez, and others
in the ALGENCAN \cite{abms07} package, with the theoretical basis for the
algorithms detailed in the book \cite{bm14}.
After the success of interior point methods for linear programming,
interior point algorithms were also developed for general nonlinear programming
including Vanderbei and Shanno's LOQO \cite{vs99}
(a merit line-search interior point method based on
a quadratic program solver also named LOQO),
Waltz and Nocedal's KNITRO \cite{WaltzNocedal2003} (an interior point
approach based on sequential quadratic programming and trust regions
\cite{Byrd06, bhn99}),
and Biegler and W\"{a}chter's IPOPT \cite{Wachter06} which employs
an interior point method with a filter line-search.
IPOPT has been adopted as a COIN-OR project
(Computational Infrastructure for Operations Research).

Performance data often shows IPOPT is among the best performing NLP solvers.
For example, in \cite{bm20} Birgin and Mart\'{\i}nez
provide performance data showing that
ALGENCAN and IPOPT are competitive with each other.
Although ALGENCAN was more robust in the experiments,
IPOPT had slightly better CPU time performance on a set of 688 problems
from CUTEst (Constrained and Unconstrained Testing Environment
with Safe Threads \cite{GouldOrbanToint15})
where both codes found equivalent solutions.
In Chapter~21 of \cite{Andrei17}, Andrei observed that for a set of 93
test problems, KNITRO and IPOPT had similar performance in terms of
number of iterations (\cite[Fig.~21.1]{Andrei17}),
while KNITRO had somewhat better performance in terms of CPU time
(\cite[Fig.~21.2]{Andrei17}).
In \cite{Wachter06}, Biegler and W\"{a}chter found that an early version
of IPOPT had better performance than an early version of KNITRO on a test set
consisting of 979 CUTE problems.
In \cite{WanBiegler} Wan and Biegler observed that in a set of 227
problems from CUTEst, KNITRO and IPOPT were very similar in performance.
KNITRO was faster than IPOPT before incorporating
regularization techniques, but slower after using regularization.
\section{Overview of PASA}
\label{overview}
%
As discussed in the introduction, PASA has two phases:
gradient projection iterations over the entire polyhedron in phase one and
projected (conjugate) gradient iterations in phase two
to optimize over faces of the polyhedron.
To choose between the two phases, we compare the violation of the
local optimality conditions for the global problem (\ref{P}) to the violation
in the local optimality conditions on the current face of the polyhedron
(the local problem).
An estimate of the violation in the optimality conditions for the
global problem is given by
\[
E (\m{x}) = \| P_\Omega (\m{x} - \nabla f (\m{x})) - \m{x} \|,
\]
where $P_\Omega$ denotes the Euclidean norm projection given by
\begin{equation}\label{proj}
P_\Omega (\m{x}) := \arg \min_{\m{y}}
\{ \|\m{x} - \m{y}\|^2 : \m{y} \in \Omega \} .
\end{equation}
Recall \cite[P7]{hz05a} that $E(\m{x}) = 0$ if and only if
$\m{x}$ is a stationary point for the global problem (\ref{P}).
After a change of variables, we obtain
\begin{equation}\label{E}
E(\m{x}) = \|\m{y}(\m{x})\| \mbox{ where }
\m{y}(\m{x}) = \arg \min_{\m{y}} \{ \|\m{y} + \nabla f (\m{x})\|:
\m{y} \in \Omega - \m{x} \} .
\end{equation}
Thus the global error at $\m{x}$ is the projection of the
negative gradient $-\nabla f (\m{x})$ onto the shifted polyhedron
$\Omega - \m{x}$.

Suppose that $\Omega$ is expressed in the form (\ref{omega}), and
for any feasible point $\m{x}$, let $\C{A}(\m{x})$ denote the active (binding)
constraints:
\[
\C{A}(\m{x}) = \{ i : (\m{Ax} - \m{b})_i = 0 \} .
\]
The active manifold at $\m{x}$ is
\[
\C{M}(\m{x}) = \m{x} + \C{N}(\m{A}_B), \quad B = \C{A}(\m{x}),
\]
where $\m{A}_B$ is the submatrix of $\m{A}$ corresponding to row indices
in $\C{A}(\m{x})$
and $\C{N}(\m{A}_B)$ is the null space of $\m{A}_B$.
The local problem corresponding to the active manifold at $\m{x}$ is
\begin{equation}\label{p}
\min f(\m{z}) \quad \mbox{subject to} \quad \m{z} \in \C{M}(\m{x}) .
\end{equation}
By (\ref{E}) with $\Omega$ replaced by $\C{M}(\m{x})$,
$\m{x}$ is a stationary point for the local problem (\ref{p})
if and only if $e(\m{x}) = 0$ where
\[
e(\m{x}) = \|\m{y}_B(\m{x})\|, \quad
\m{y}_B(\m{x}) =
\arg \min_{\m{y}} \{ \|\m{y} + \nabla f (\m{x})\|: \m{A}_B\m{y} = \m{0}\} .
\]
Thus the local error bound at $\m{x}$ is the projection of the
negative gradient $-\nabla f (\m{x})$ onto $\C{N}(\m{A}_B)$.

\begin{figure}
\begin{center}
\includegraphics[scale=.33]{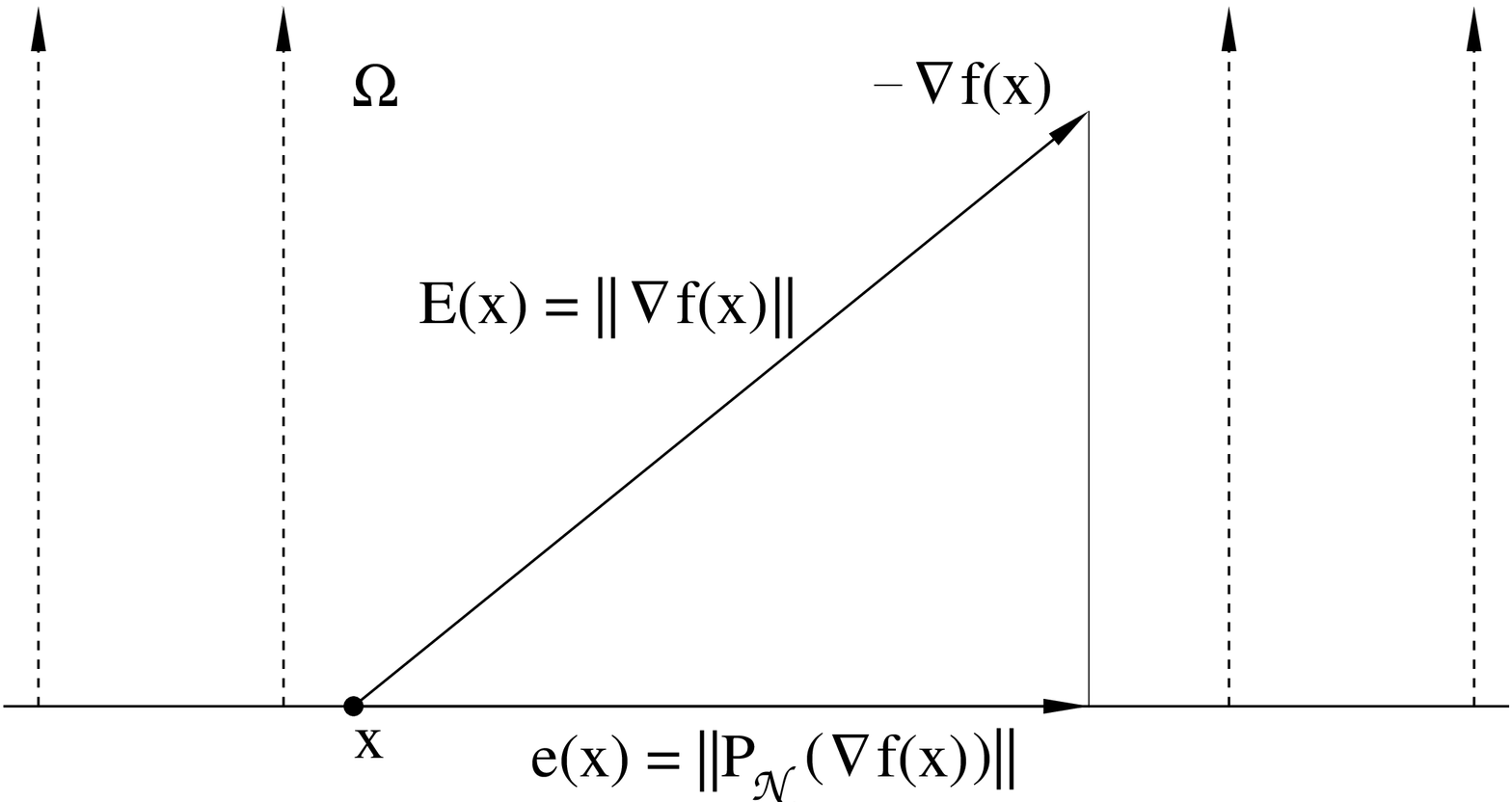}
\caption{Global error $E(x)$ versus local error $e(x)$.
\label{l_versus_g}}
\end{center}
\end{figure}

A simple illustration of the local and global errors is given in
Figure~\ref{l_versus_g}.
The set $\Omega$ is the upper half-space and the point $\m{x}$ lies on the
boundary of $\Omega$.
Since the negative gradient at $\m{x}$ points into $\Omega$,
$E(\m{x})$ is simply the norm of the negative gradient.
The active manifold is the horizontal axis, and $e(\m{x})$ is the
projection of the negative gradient onto the horizontal axis.

\renewcommand\figurename{Alg.}
\begin{figure}[h]
\begin{center}
{ \tt
\begin{tabular}{l}
\hline
\\
{\bf \textcolor{blue}{Parameters:}} \hspace*{.02in}$\theta \in (0, 1)$,
$\tau \in (0, \infty)$,
start guess $\m{x}_0 \in \mathbb{R}^n$\\
{\bf \textcolor{blue}{Initialization:}}\hspace*{.000in}
$\m{x}_1 = \C{P}_\Omega (\m{x}_0)$, $k = 1$ \\
{\bf \textcolor{blue}{Phase one:}}
\hspace*{.10in}While $E(\m{x}) > \tau$, execute phase one\\
\hspace*{1.25in}Possibly reduce $\theta$\\
\hspace*{1.25in}If $e(\m{x}_k) \ge \theta E(\m{x}_k)$, goto phase two; \\
\hspace*{1.25in}else $k \leftarrow k+1$.
\\
\hspace*{1.00in}End\\[.1in]
{\bf \textcolor{blue}{Phase two:}}
\hspace*{.10in}While $E(\m{x}) > \tau$, execute phase two\\
\hspace*{1.25in}Possibly reduce $\theta$\\
\hspace*{1.25in}If $e(\m{x}_k) < \theta E(\m{x}_k)$, goto phase one; \\
\hspace*{1.25in}else $k \leftarrow k+1$.  \\
\hspace*{1.00in}End\\
\hline
\end{tabular}
}
\end{center}
\caption{Sketch of Polyhedral Active Set Algorithm (PASA).
\label{pasa}}
\end{figure}
\renewcommand\figurename{Fig.}

In implementing PASA, we choose a parameter $\theta \in (0,  1)$,
and then operate in either phase one or phase two as indicated in
Algorithm~\ref{pasa}.
As seen in Algorithm~\ref{pasa}, the branching between the two phases
of PASA is based on a comparison between $E$ and $e$ at the current
iterate $\m{x}_k$.
It is shown in \cite{HagerZhang16} that
\[
\liminf\limits_{k \to\infty} E(\m{x}_k) =0,
\]
whenever the algorithm in phase two satisfies the following conditions:
\begin{enumerate}
\item [P1.]
For each $k$, $\m{x}_{k} \in \Omega$ and $f(\m{x}_{k+1}) \le f(\m{x}_k)$.
\item [P2.]
For each $k$, $\C{A}(\m{x}_{k}) \subset \C{A}(\m{x}_{k+1})$.
\item [P3.]
If $\C{A}(\m{x}_{j+1}) = \C{A}(\m{x}_j)$
for $j \ge k$, then $\liminf\limits_{j\to\infty} e(\m{x}_j) = 0$.
\end{enumerate}
Moreover, under either a nondegeneracy or strong second-order
sufficient optimality with linear independence of the active
constraint gradients,
the iterates of PASA are only generated by phase two
when $k$ is sufficiently large.
To achieve this stronger property, an adjustment was made in \cite{HagerZhang16}
of the form $\theta \leftarrow \mu \theta$, $\mu \in (0, 1)$,
whenever the ``undecided index set was empty'' in phase one.
The undecided indices corresponded to constraints for which the associated
multipliers were either sufficiently positive or the constraint was
sufficiently active.
Determining whether the undecided index set was empty required an
additional projection which detracted from the efficiency of PASA.
Note that any update to $\theta$ in phase one which drives it to zero
guarantees that the iterates are generated by phase two asymptotically
under the condition given above.
Since phase one often branches to phase two after a single iteration,
a practical approach for driving $\theta$ to zero when too much time is
spent in phase one is to
decrease $\theta$ whenever more than one iteration is performed
in phase one.

The rules for branching between phases one and two are based on the following
considerations.
First note that phase one, by itself, is typically globally convergent
(see \cite{hz05a}).
The convergence rate, however, is at best linear.
The purpose of phase two is to improve efficiency by using a superlinearly
convergent algorithm to find an optimum over the manifold defined by the
active constraints.
Nonetheless, efficiency is lost if the optimum over the manifold is computed
with too much precision.
Since the default value for $\theta$ is 0.01, we would branch from phase two
to phase one if the local error $e(\m{x}_k)$ is less than 0.01 times the
global error $E(\m{x}_k)$.
In phase one, typically only one iteration is performed,
some constraints that were active become inactive, and the iterate
moves to a new manifold.
In the previous iteration, the inequality $e(\m{x}_k) < \theta E(\m{x}_k)$
was satisfied.
After performing the gradient projection step of phase one and moving
to a new active manifold, we usually find that 
$e(\m{x}_k) > \theta E(\m{x}_k)$, so PASA branches back to phase two and
begins to explore a new active manifold.
%
\section{Phase One}
\label{phase_one}
\begin{figure}
\begin{center}
\includegraphics[scale=.4]{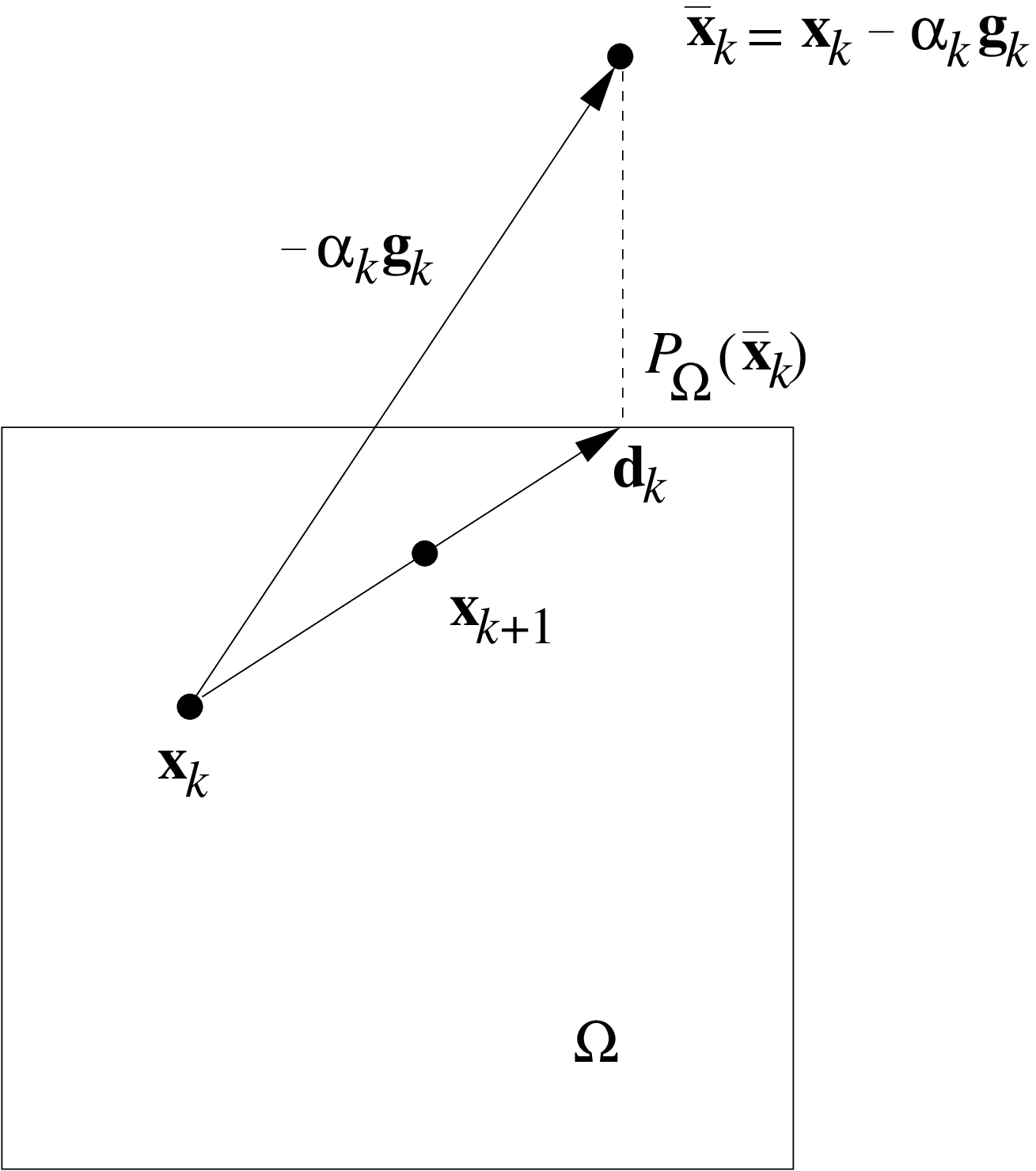}
\caption{Sketch of the gradient projection algorithm.
\label{gradproj}}
\end{center}
\end{figure}

%
The version of the gradient projection algorithm that we utilize
is depicted in Figure~\ref{gradproj}.
At the current iterate $\m{x}_k$, a step of length $\alpha_k$ is
taken along the negative gradient $-\m{g}_k = -\nabla f (\m{x}_k)$
to reach a point $\bar{\m{x}}_k$,
whose projection onto the polyhedron $\Omega$ is
$P_\Omega (\bar{\m{x}}_k)$.
A line search is performed along the search direction
$\m{d}_k = P_\Omega (\bar{\m{x}}_k) - \m{x}_k$ to obtain the
next iterate $\m{x}_{k+1}$.
Our choice for the stepsize $\alpha_k$ is based on the BB formula \cite{bb88}
with a cyclic implementation which is explained in \cite{dhsz05}.
In a cyclic implementation, the stepsize is kept fixed in some iterations
(that is, $\alpha_{k+1} = \alpha_k$), while in other iterations,
it is given by the BB formula.
As shown in \cite{dhsz05}, a cyclic implementation can lead to better
performance.

Our version of the gradient projection algorithm is similar to the version in
SPG \cite{bmr01} since the line search is along the line segment connecting
$\m{x}_k$ to the projected point $P_\Omega (\bar{\m{x}}_k)$.
There is another version of the gradient projection method in which
$\m{x}_{k+1} = P_\Omega (\m{x}_k - s_k \nabla f (\m{x}_k))$
where the stepsize $s_k$ is chosen to satisfy both a descent condition
and a condition to ensure the $s_k$ is not too small;
for example, see \cite{bm88, bm94, bmt90}
and the references therein.
For this scheme,
the active constraints at a minimizer can be identified
in a finite number of iterations under suitable assumptions, however,
more than one projection may be needed to determine an acceptable step.
Although the gradient projection algorithm implemented in PASA may not
identify the active constraints, we show in \cite[Lemma~6.2]{HagerZhang16}
that the violation in the active constraints is on the order of the
squared error in the iterate $\m{x}_k$.

The line search implemented in PASA,
shown in Algorithm~\ref{phase1}, is of Armijo type \cite{armijo}.
If $f_k^r = f(\m{x}_k)$, then this is
an ordinary Armijo line search restricted to the feasible set $\Omega$.
However, the implemented line search is nonmonotone, such as in \cite{gll86},
but with a more sophisticated choice of $f_k^r$,
based on the procedure given in the appendix of \cite{hz05a}.
Also, to avoid potential breakdown of the line search in a neighborhood
of an optimum, an approximate (but more accurate) line search
is used near a local minimizer;
see the approximate Wolfe line search in \cite{hz03}.
\renewcommand\figurename{Alg.}
\begin{figure}[t]
\begin{center}
{ \tt
\begin{tabular}{l}
\hline
\\
{\bf \textcolor{blue}{Parameters:}}
$\delta$ and $\eta \in (0, 1)$, $\alpha_k \in (0, \infty)$\\[.05in]
While $E(\m{x}_k) >$ max \{$\tau$, $e(\m{x}_k)/\theta$\}\\[.05in]
\hspace*{.3in}1.~$\m{d}_k =$
$\C{P}_\Omega (\m{x}_k - \alpha_k \m{g}_k) - \m{x}_k$ \\[.05in]
\hspace*{.3in}2.~$s_k = \eta^j$ where $j \ge 0$ is smallest integer
such that \\[.05in]
\hspace*{.6in}$f(\m{x}_k + s_k \m{d}_k) \le f_k^r + s_k \delta
\nabla f(\m{x}_k)\m{d}_k$\\[.05in]
\hspace*{.3in}3.~$\m{x}_{k+1} = \m{x}_k + s_k \m{d}_k$
and $k \leftarrow k + 1$
\\[.05in]
End\\
\hline
\end{tabular}
}
\end{center}
\caption{Phase one (gradient projection algorithm). \label{phase1}}
\end{figure}
\renewcommand\figurename{Fig.}

%
\section{Phase Two}
\label{phase_two}
%
According to the theory developed in \cite{HagerZhang16},
any algorithm with the properties (P1)--(P3) can be used in phase two.
The current gradient-based implementation of PASA combines an
{\it active set} gradient projection algorithm with a {\it projected}
version of the conjugate gradient method.
Let us define the set
\[
\Omega_k = \{ \m{x} \in \Omega : (\m{Ax}-\m{b})_i = 0 \mbox{ for all }
i \in \C{A}(\m{x}_k) \} .
\]
By an active set gradient projection algorithm (A-GP), we mean that
\[
\m{x}_{k+1} = \m{x}_k + s_k \m{d}_k, \quad \mbox{where }
\m{d}_k = \C{P}_{\Omega_k} (\m{x}_k - \alpha_k \m{g}_k) - \m{x}_k ,
\quad 0 < s_k \le 1.
\]
This is the gradient projection step of Algorithm~\ref{phase1} except that
$\Omega$ is replaced by $\Omega_k$.
Since the constraints that are active at $\m{x}_k$ are also active at
$\C{P}_{\Omega_k} (\m{x}_k - \alpha_k \m{g}_k)$,
all the constraints active at $\m{x}_k$ are also active at $\m{x}_{k+1}$.
Potentially, when $s_k = 1$, additional constraints could be active at
$\m{x}_{k+1}$.
As long as $\C{A}(\m{x}_{k})$ is strictly contained in
$\C{A}(\m{x}_{k+1})$ and
phase two does not reach one of the termination conditions
in Algorithm~\ref{pasa}, A-GP continues to operate.
At any iterate $\m{x}_k$ where
$\C{A}(\m{x}_k) = \C{A}(\m{x}_{k-1})$,
we switch in phase two from A-GP to a conjugate gradient scheme.
This switch is done to exploit the faster convergence of the conjugate
gradient method when compared to gradient descent.
When the active set is growing, it is pointless to switch to conjugate
gradients since CG needs to be restarted whenever a new constraint becomes
active.
Hence, the switch to CG is not attempted until
$\C{A}(\m{x}_k) = \C{A}(\m{x}_{k-1})$ in A-GP.

Our implementation of conjugate gradients in phase two is now explained.
Let $\m{A}_k$ be the submatrix of $\m{A}$ associated with the
active constraint gradients at $\m{x}_k$, and let $\m{b}_k$ denote
the associate right side of the constraint.
During phase two, the constraint $\m{A}_k \m{x} \le \m{b}_k$ is enforced
as an equality, while the remaining constraints are strict inequalities
at $\m{x}_k$.
Since $\m{A}_k \m{x}_k = \m{b}_k$, the change of variables
$\m{x} = \m{x}_k + \m{z}$ yields the equation $\m{A}_k \m{z} = \m{0}$.
After this change of variables, the optimization problem in phase two
is written
\begin{equation}\label{Ak}
\min f(\m{x}_k + \m{z}) \quad \mbox{subject to} \quad
\m{A}_k\m{z} = \m{0}, \quad \m{z} \in \Omega - \m{x}_k.
\end{equation}

If $\m{P}_k \in \mathbb{R}^{n \times n}$ denotes the
orthogonal projection onto the null space $\C{N}(\m{A}_k)$,
then the change of variable $\m{z} = \m{P}_k \m{y}$ in (\ref{Ak})
yields the locally unconstrained problem
\begin{equation}\label{Uk}
\min f(\m{x}_k + \m{P}_k\m{y}) \quad \mbox{subject to} \quad
\m{x}_k + \m{P}_k \m{y} \in \Omega.
\end{equation}
This problem is locally unconstrained since $\m{y} = \m{0}$ lies in the
interior of the set $\Omega - \m{x}_k$.
Hence, the conjugate gradient algorithm can be applied to (\ref{Uk})
generating iterates $\m{y}_j$, $j \ge 0$, starting from $\m{y}_0 = \m{0}$.

Each conjugate gradient iteration involves a search direction $\m{d}_j$
and a line search along $\m{d}_j$.
If $\alpha$ is the stepsize along $\m{d}_j$, the line search
enforces the constraint
\begin{equation}\label{Ik}
\m{x}_k + \m{P}_k(\m{y}_j + \alpha \m{d}_j) \in \Omega.
\end{equation}
Assuming $\m{P}_k \m{y}_j$ lies in the interior of the feasible set
$\Omega - \m{x}_k$, we let $\alpha_{\max}$ denote the largest $\alpha$
such that the inclusion (\ref{Ik}) holds.
With this notation, the conjugate gradient line search focuses on the
optimization problem
\[
\min \;\; f(\m{x}_k + \m{P}_k(\m{y}_j + \alpha \m{d}_j))
\quad \mbox{subject to} \quad
0 \le \alpha \le \alpha_{\max} .
\]
If the solution of this problem is $\alpha = \alpha_{\max}$, then
one or more constraints are activated and we return to A-GP.
Otherwise, the projected conjugate gradient iteration continues.
Phase two is summarized in Algorithm~\ref{phase2}.

\renewcommand\figurename{Alg.}
\begin{figure}
\begin{center}
{ \tt
\begin{tabular}{l}
\hline
\\
While $e(\m{x}_k)/\theta \ge E(\m{x}_k) > \tau$ \\[.05in]
\hspace*{.3in}1.~Perform A-GP until
$\C{A}(\m{x}_k) = \C{A}(\m{x}_{k-1})$, then branch to step 2.\\[.05in]
\hspace*{.3in}2.~Apply limited memory CG\_DESCENT to (\ref{Uk});
branch to\\[.05in]
\hspace*{.52in}step 1 when reaching $\m{x}_{k+1} =$
$\m{x}_k + \m{P}_k(\m{y}_j + \alpha_{\max} \m{d}_j)$
on the\\[.05in]
\hspace*{.52in}boundary of $\Omega$.
\\[.05in]
End\\
\hline
\end{tabular}
}
\end{center}
\caption{Phase two (A-GP and CG\_DESCENT). \label{phase2}}
\end{figure}
\renewcommand\figurename{Fig.}

A stability issue arises when applying the conjugate gradient method
to the unconstrained problem (\ref{Uk}).
Theoretically, the projection can be expressed as
\[
\m{P}_k = \m{I} - \m{A}_k\tr(\m{A}_k \m{A}_k\tr)^{-1}\m{A}_k,
\]
where the inverse should be replaced by a pseudoinverse when the
rows of $\m{A}_k$ are linearly dependent, and $\m{P}_k$ is a positive
semidefinite matrix.
The routine PPROJ, used in phase one computes the factorization
\[
\m{A}_k \m{A}_k\tr + \sigma \m{I} = \m{LDL}\tr ,
\]
where $\sigma > 0$ is relatively small, $\m{L}$ is lower triangular with
ones on the diagonal, and $\m{D}$ is diagonal.
This leads us to replace $\m{P}_k$ in (\ref{Uk}) by its approximation
\[
\tilde{\m{P}}_k =
\m{I} - \m{A}_k\tr(\m{A}_k\m{A}_k\tr + \sigma \m{I})^{-1}\m{A}_k =
\m{I} - \m{A}_k\tr(\m{LDL})^{-\sf T}\m{A}_k,
\]
a positive definite matrix.

Since our goal is to optimize the objective over
$\m{x} = \m{x}_k + \m{P}_k\m{y}$,
it is convenient to formulate the algorithm
for solving (\ref{Uk}) in terms of $\m{x}$ rather than in terms of $\m{y}$.
In particular, when the iterates are given by the
CG\_DESCENT family parameterized by $\eta > 1/4$,
the search directions are (see \cite[Sect. 2]{HagerZhang13c})
\begin{equation}\label{d}
\m{d}_0 = -\tilde{\m{P}}_k^2 \m{g}_0, \quad
\m{d}_{k+1} = -\tilde{\m{P}}_k^2 \m{g}_{k+1} + \beta_k \m{d}_k \; \;
\mbox{ for } k \ge 0,
\end{equation}
where $\m{g}_k = \nabla f(\m{x}_k)$ and
\begin{equation}\label{beta}
\beta_k = \frac{\m{y}_k\tr  \tilde{\m{P}}_k^2 \m{g}_{k+1}}{\m{d}_k\tr \m{y}_k}
- \eta \frac{\|\tilde{\m{P}}_k\m{y}_k\|^2}{\m{d}_k\tr\m{y}_k}
\frac{\m{d}_k\tr\m{g}_{k+1}}{\m{d}_k\tr\m{y}_k}, \quad
\m{y}_k = \m{g}_{k+1} - \m{g}_k .
\end{equation}

Even though $\tilde{\m{P}}_k$ is an approximation to $\m{P}_k$ and
$\m{P}_k^2 = \m{P}_k$, the replacement of $\tilde{\m{P}}_k^2$ by
$\tilde{\m{P}}_k$ leads to very poor performance.
Moreover, when we compute the search directions by the formula
(\ref{d}), the iterates quickly lose feasibility.
The reason is that the component of the error in $\m{d}_k$ pointing
out of $\C{N}(\m{A}_k)$ is added into $\m{d}_{k+1}$ in (\ref{d}),
and these errors in the search direction can accumulate.
The following iteration is equivalent to (\ref{d}) and numerically stable
since $\m{d}_{k+1}$ is the product of an intermediate
vector $\m{D}_{k+1}$ with $\tilde{\m{P}}_k$ which removes error components
orthogonal to the null space of $\m{A}_k$:
\begin{equation}\label{D}
\m{D}_0 = -\tilde{\m{P}}_k \m{g}_0, \quad
\m{D}_{k+1} = -\tilde{\m{P}}_k \m{g}_{k+1} + \beta_k \m{D}_k, \quad
\m{d}_{k+1} = \tilde{\m{P}}_k \m{D}_{k+1}  \mbox{ for } k \ge 0.
\end{equation}
%
\section{Results}
\label{results}
%
We compare the performance of PASA Version 2.0.0 to the performance of
IPOPT Version 3.14.5 using the CUTEst platform \cite{GouldOrbanToint15},
and polyhedral constrained optimization problems from CUTEst along with the
Maros/Meszaros quadratic programming test set.
IPOPT can operate in a gradient-based mode, where 
the Hessian of the Lagrangian in the KKT system is approximated by a limited
memory quasi-Newton method (L-BFGS),
and a Hessian-based mode when both
the gradient and Hessian of the objective and the constraints are
provided, and a direct solver is used for the linear systems.
We installed both of the recommended linear solvers:
MUMPS 5.4.1 and the HSL software, which includes MA57, Version 3.11.0.
When IPOPT was run, it always chose the MUMPS linear solver.
In comparisons between the gradient and Hessian-based IPOPT,
the Hessian-based version performed much better.
Hence, our comparisons are with Hessian-based IPOPT.
Note that comparisons between the Hessian-based PASA
(currently under development) and the gradient-based PASA also
indicate that Hessian-based PASA is superior to gradient-based PASA.

In selecting the problem set for the numerical experiments, 42 of the
QPs from CUTEst were excluded.
The names of these problems begin with the letter A followed by either
0 or 2 or 5.
There were two issues with this subset of the CUTEst test set.
First, in many cases, the starting point is essentially a stationary point,
and gradient-based PASA immediately terminates.
Second, these problems have between 15,000 and 20,000 linear constraints,
and a small number of dense columns.
Due to the dense columns,
the matrix $\m{AA}\tr$ is dense with dimension between 15,000 and 20,000.
To handle these problems efficiently, the dense columns need to be removed and
processed using a Woodbury update \cite{Hager89}.
We have not yet had time to incorporate Woodbury updates in PASA.
Moreover, if the Woodbury updates were incorporated in the code,
termination may occur at the starting point, and the problem would be
excluded by the rules given in the next paragraph.
Note that the Hessian-based PASA should handle these problems without
difficulty since the KKT system has a sparse factorization.

After these exclusions, we start with 655 problems which we tried to solve
to the accuracy tolerance 1.e$-$6.
If the objective values computed by each solver agreed to
4 significant digits, then we accepted the problem.
If 4 digit agreement was not achieved, then we examined the computed solutions.
If the solvers were converging to different solutions, then we removed
the problem from the test set; in other words, we focused on problems where both
solvers started from the same initial guess and reached the same solution.
There were 75 problems where the solvers converged to different solutions.
In 38 cases, the solution computed by IPOPT had a better objective value,
and in 37 cases, the solution computed by PASA had a better objective value.
Note that among the 38 cases where IPOPT had a better objective value,
it was observed that in a number of these cases, the starting guess was
essentially a stationary point, and PASA stopped immediately,
while the Hessian-based IPOPT did not stop at the starting point.
Both solvers, however, are only guaranteed to converge to a stationary point.

After pruning the 75 problems where the solvers converged to different
solutions, there were 580 remaining test problems.
If the objective values disagreed by more than 4 significant digits but
the solvers were converging to the same solution, we then adjusted the
accuracy tolerance of the less accurate solver so as to achieve
comparable accuracy to that of the more accurate solver.
In these cases where one solver was more accurate than the other,
we found that the PASA estimate $E(\m{x})$ for the solution tolerance
resulted in a more accurate objective value in most cases.
When the accuracy tolerance of IPOPT was adjusted to match the accuracy of
PASA, often just one or two more iterations were needed.
When a solver was unable to achieve the accuracy tolerance 1.e$-$6
for a problem, its computing time was set to $\infty$.


The performance of the gradient-based PASA and Hessian-based IPOPT are
compared using wall time.
Note that Hessian-based algorithms such as either IPOPT or the Hessian-based
PASA, typically require fewer iterations and evaluations
(function and gradient) when compared to the gradient-based PASA.
The trade-offs between convergence rate of an algorithm
and evaluation time are not studied in this paper,
instead we focus on wall time.
The run data is available at:
\smallskip
\begin{center}
https://people.clas.ufl.edu/hager/files/IPOPTresults.txt
https://people.clas.ufl.edu/hager/files/PASAresults.txt
\end{center}
\smallskip

Both solvers exploit multiple processors when matrices are factored and
linear systems are solved.
The software was run on a Lenovo ThinkPad with 8 Intel
i7-865U CPUs operating at 1.90GHz (4 cpu cores) with 8,192~KB cache and
16~GB memory.
The operating system was Ubuntu Linux
with Intel's MKL (Math Kernel Library) BLAS.
PASA used the timer {\it gettimeofday} with microsecond accuracy, while
IPOPT appears to use the timer {\it ftime} with millisecond accuracy
(embedded inside the function IpCoinGetTimeOfDay).
\begin{figure}
\begin{center}
\includegraphics[scale=.28]{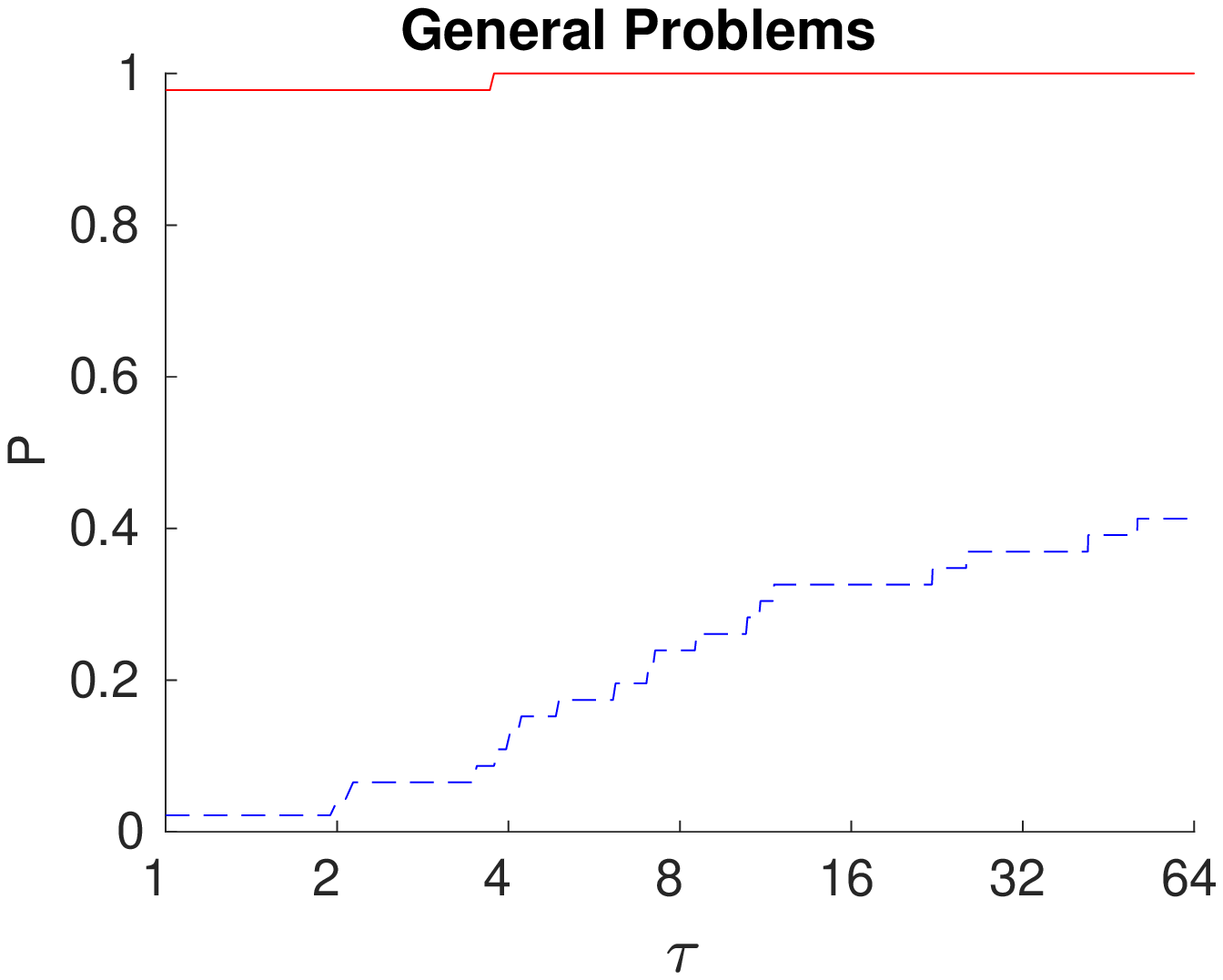}
\includegraphics[scale=.28]{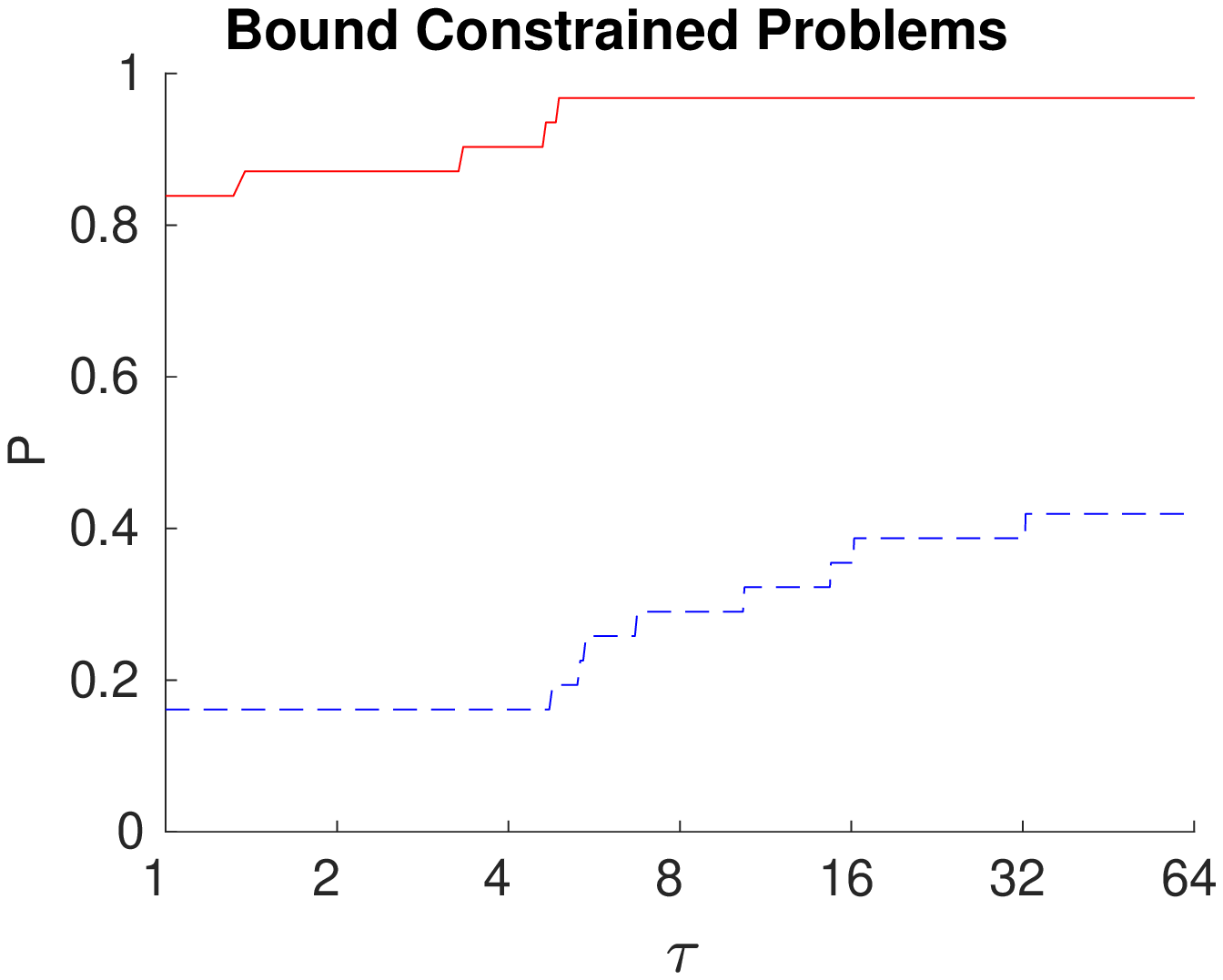}
\includegraphics[scale=.28]{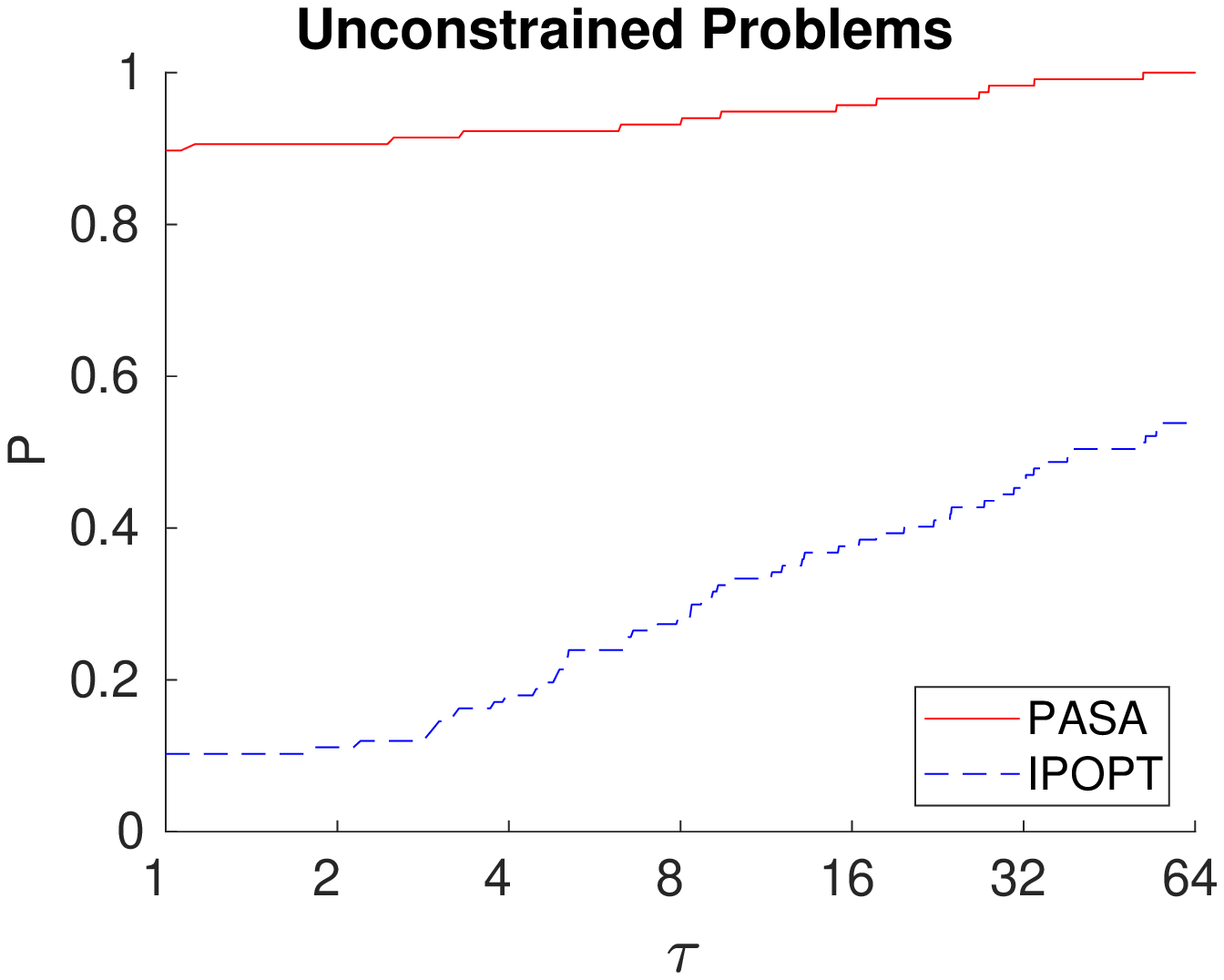}
\caption{Wall time performance profiles for general, bound,
and unconstrained programs
\label{general}}
\end{center}
\end{figure}
\begin{figure}
\begin{center}
\includegraphics[scale=.28]{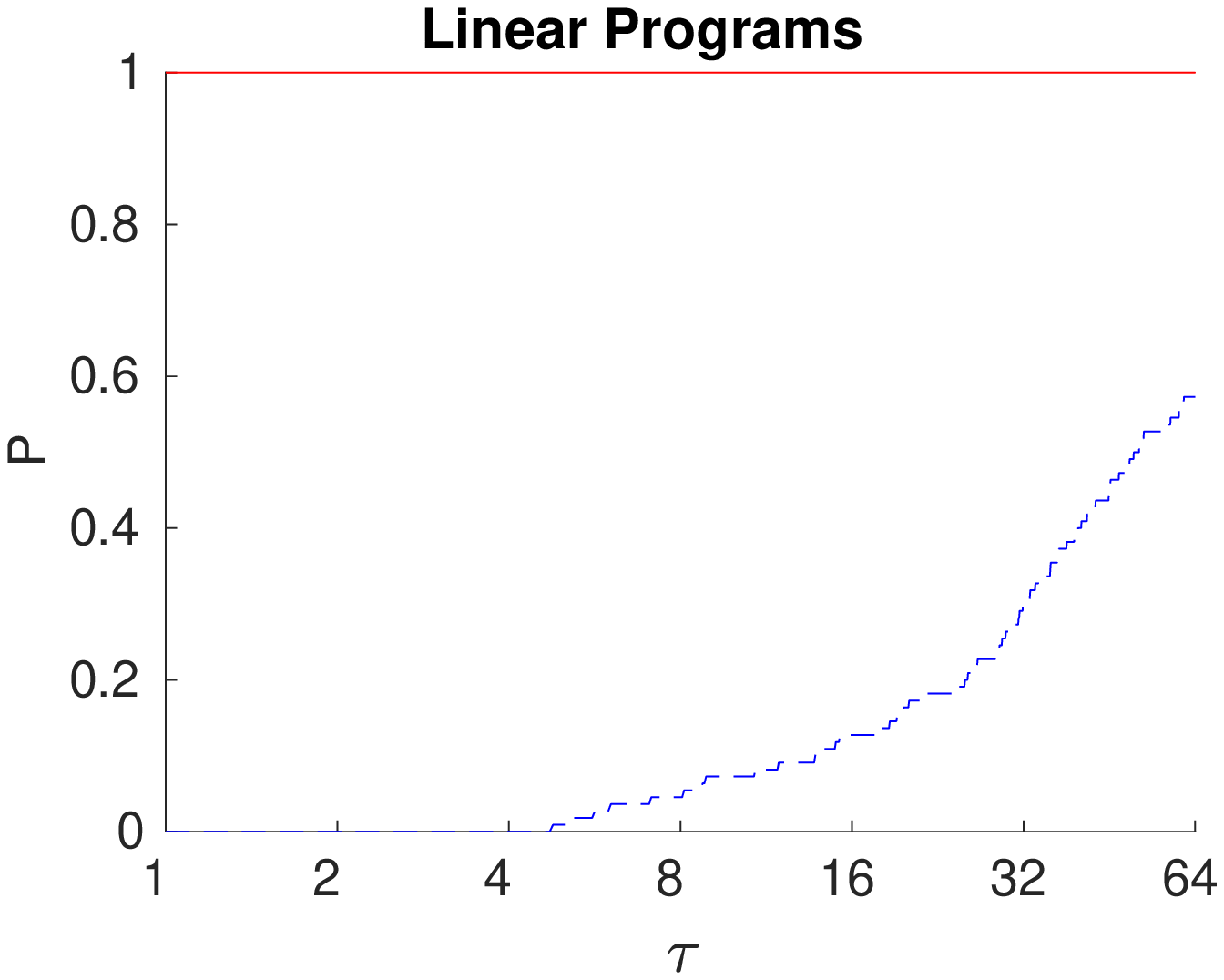}
\includegraphics[scale=.28]{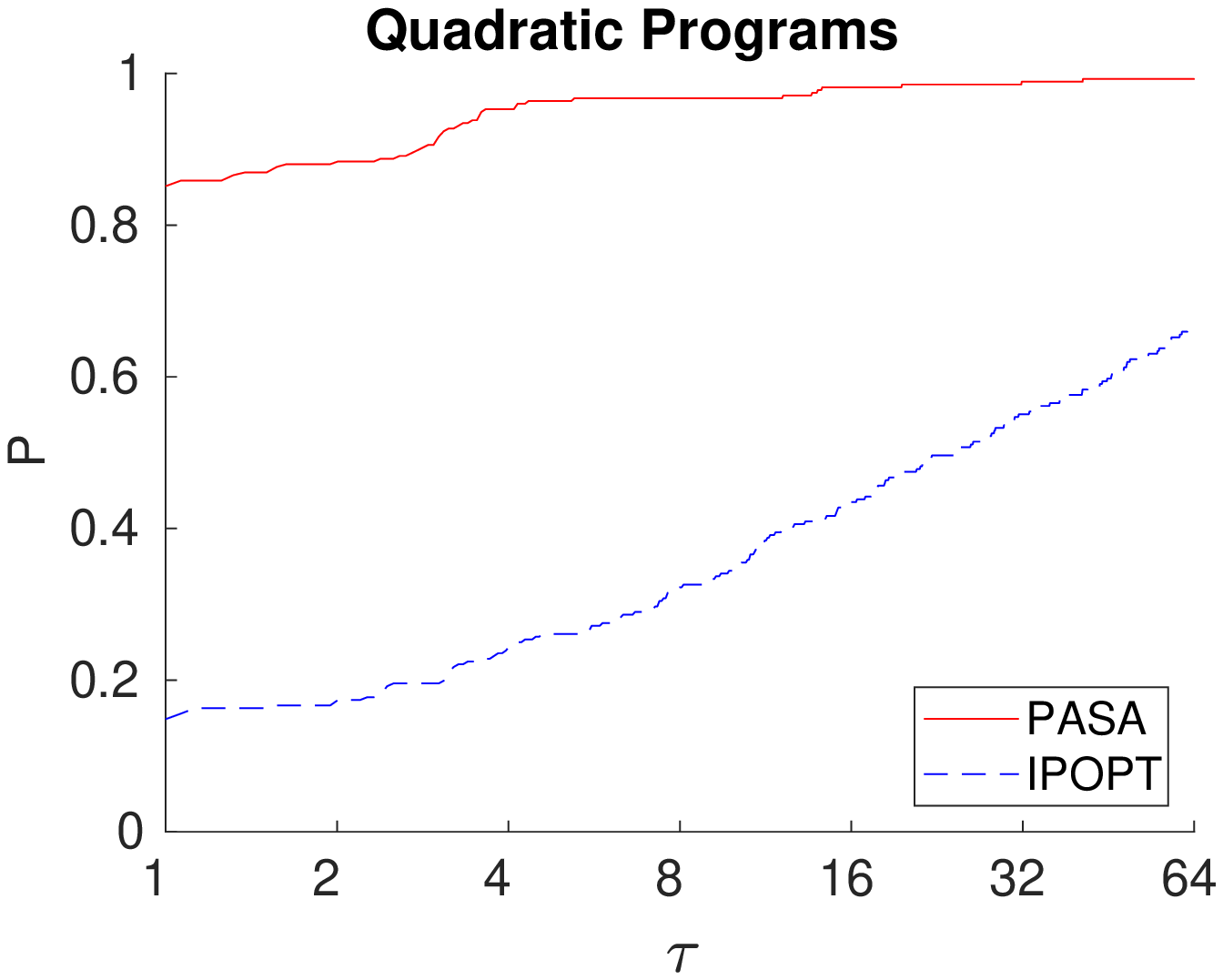}
\includegraphics[scale=.28]{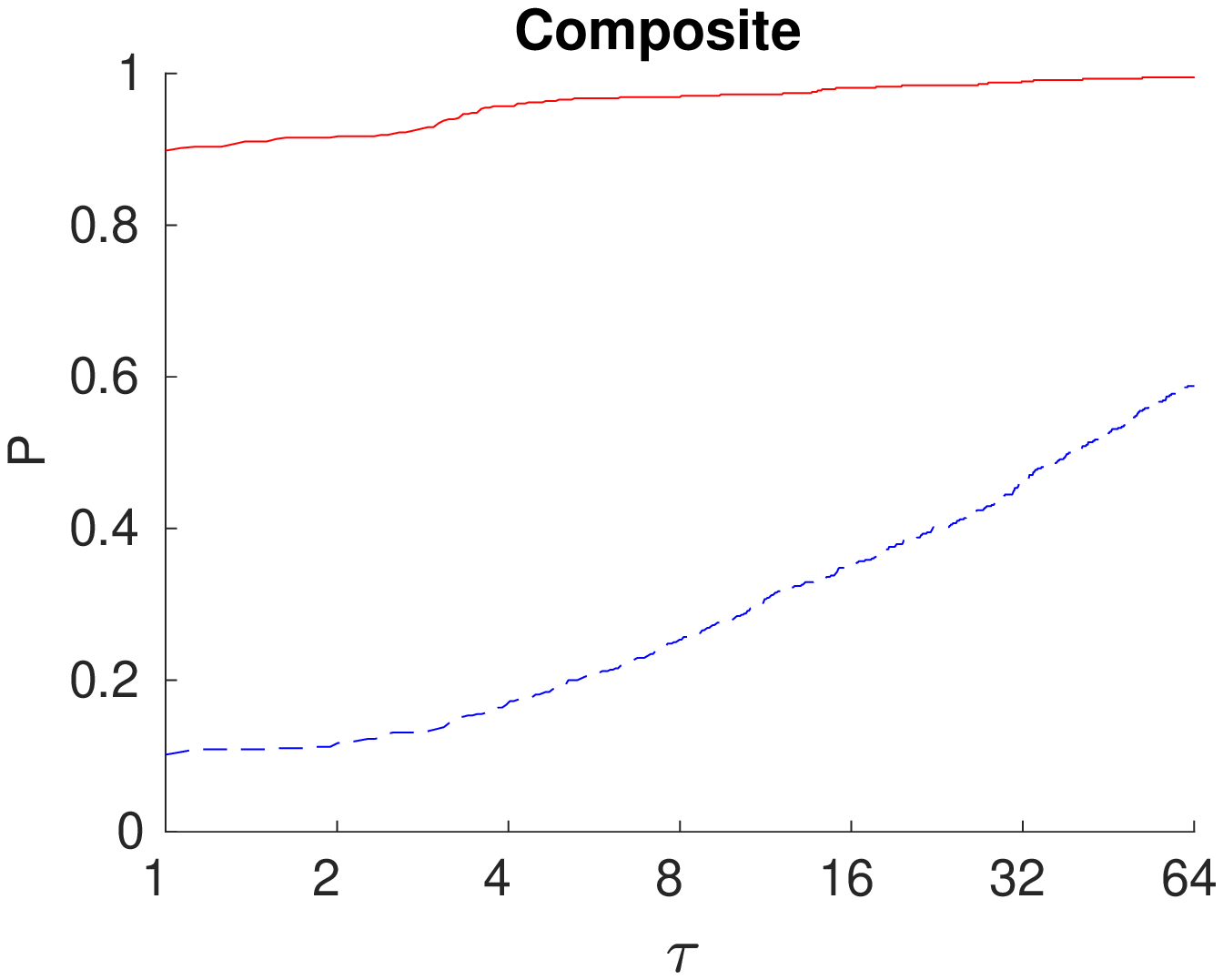}
\caption{Wall time performance profiles for linear, quadratic,
and composite programs
\label{linear}}
\end{center}
\end{figure}

Figures~\ref{general}--\ref{linear} plot the wall time performance
profiles for the two codes.
The vertical axis gives the fraction P of problems for which
any given method is within a factor $\tau$ (horizontal axis) of the best time.
The top curve is the method that solved the most problems in a time
that was within a factor $\tau$ of the best time.
The percentage of the test problems for which a method is fastest is
given on the left axis of the plot.
The right side of the plot gives the percentage of the test problems
that were successfully solved by each of the methods.
In essence, the right side is a measure of an algorithm's robustness.

In preparing the plots, the problems in the test set were partitioned
into 5 groups:
linear and quadratic programs denote polyhedral constrained problem
for which the objective is linear or quadratic respectively.
General problems have nonquadratic nonlinear objectives
with additional linear and possibly bound constraints.
Bound constrained problems have nonquadratic nonlinear objectives
and only bound constraints.
Unconstrained problems have nonquadratic nonlinear objectives without
constraints.
The composite problems are the union of all 5 groups.
Based on the plots, gradient-based PASA performed relatively well on this
collection of test problems where function and gradient evaluations are
relatively cheap; the cost of the linear algebra in IPOPT for solving
the linear systems of equations outweighed the savings associated
with a lower number of evaluations.

When a problem is unconstrained, $E(\m{x}) = e(\m{x})$ and since $\theta < 1$,
PASA immediately branches to phase two, where it remains until the
convergence tolerance is satisfied.
Hence, the performance on unconstrained problems essentially reflects
the performance of limited memory CG\_DESCENT \cite{HagerZhang13c}.
For bound constrained problems, PASA's local and global error estimators
$e(\m{x})$ and $E(\m{x})$ reduce to the same estimators that were used in
the algorithm \cite{hz05a} for bound constrained problems.
Hence, the performance on bound constrained problem essentially reflects the
performance of the active set algorithm \cite{hz05a}.
Linear programs are solved in PASA by a series of gradient projection steps,
where the stepsize choice is crucial.
The performance corresponds to the first-order algorithm in \cite{DavisHager08}.
Details will be provided in a separate paper.


\section{Conclusion}
\label{conclusion}
A gradient-based implementation of the Polyhedral Active Set Algorithm (PASA)
was presented.
The algorithm was composed of two phases, the gradient projection algorithm
was used in phase one, while phase two optimized the objective over
faces of the polyhedron.
Branching between phases was determined by the relationship between
local and a global error estimators $e$ and $E$ respectively.
At a feasible point $\m{x}$ for the polyhedron, we branch from phase one
to phase two when $e(\m{x}) \ge \theta E(\m{x})$, where $\theta \in (0,1)$
is a given parameter;
we branch from phase two to phase one when $e(\m{x}) < \theta E(\m{x})$.
With suitable adjustments to $\theta$,
the iterates perform phase two asymptotically.
It was found that PASA had significantly better wall time performance
when compared to IPOPT using a collection of 580 test problems taken from
both CUTEst and the Maros/Meszaros quadratic programming test set.
Even though Hessian-based IPOPT used significantly fewer evaluations of
the objective and gradient when compared to gradient-based PASA,
the time for the linear algebra in IPOPT outweighed the savings
derived from the fewer evaluations in the test set.
\section{Acknowledgements}
The assistance of Nicholas Gould and Dominique Orban in configuring CUTEst
to enable its operation with PASA was greatly appreciated.
An initial draft of the PASA MATLAB interface by James Diffenderfer is
gratefully acknowledged.
%
\bibliographystyle{siam}

\end{document}